\magnification1200
\input amssym.def 
\input amssym.tex 
\def\SetAuthorHead#1{}
\def\SetTitleHead#1{}
\def\NoindentAfter{\everypar={\setbox0=\lastbox\everypar={}}}
\def\H#1\par#2\par{{\baselineskip=15pt\parindent=0pt\parskip=0pt
 \leftskip= 0pt plus.2\hsize\rightskip=0pt plus.2\hsize
 \bf#1\unskip\break\vskip 4pt\rm#2\unskip\break\hrule
 \vskip40pt plus4pt minus4pt}\NoindentAfter}
\def\HH#1\par{{\bigbreak\noindent\bf#1\medskip}\NoindentAfter}
\def\HHH#1\par{{\bigbreak\noindent\bf#1\unskip.\kern.4em}}
\def\th#1\par{\medbreak\noindent{\bf#1\unskip.\kern.4em}\it}
\def\endth{\medbreak\rm}
\def\pf#1\par{\medbreak\noindent{\it#1\unskip.\kern.4em}}
\def\df#1\par{\medbreak\noindent{\it#1\unskip.\kern.4em}}
\def\enddf{\medbreak}
\let\rk\df\let\endrk\enddf
\let\Roster\bgroup\def\endRoster{\egroup\par}
\def\\{}\def\text#1{\hbox{\rm #1}}
\def\mop#1{\mathop{\rm\vphantom{x}#1}\nolimits}
\def\MaxReferenceTag#1{}
\def\qedbox{\vrule width2mm height2mm\hglue1mm\relax}
\def\qed{\ifmmode\qedbox\else\hglue5mm\unskip\hfill\qedbox\medbreak\fi\rm}
\let\SOverline\overline

\let\Item\item\let\ItemItem\itemitem
\def\cite#1{{\bf[#1]}}
\def\Em#1{{\it #1\/}}
\def\Bib#1\par{\bigbreak\bgroup\centerline{#1}\medbreak\parindent30pt
 \parskip2pt\frenchspacing\par}
\def\endBib{\par\egroup}
\newdimen\Overhang
\def\rf#1{\par\noindent\hangafter1\hangindent=\parindent
     \setbox0=\hbox{[#1]}\Overhang\wd0\advance\Overhang.4em\relax
     \ifdim\Overhang>\hangindent\else\Overhang\hangindent\fi
     \hbox to \Overhang{\box0\hss}\ignorespaces}

\def\bbQ{{\Bbb Q}}
\def\bbZ{{\Bbb Z}}
\newcount\notenumber\notenumber0
\def\note#1{\advance\notenumber1\footnote{$^{\the\notenumber)}$}{#1}}
\def\Coordinates{\bigbreak\bgroup\parindent=0pt\obeylines}
\def\endCoordinates{\egroup}
\overfullrule=0pt
\let\ProofingOff\relax
\let\skipaline\smallbreak

\ProofingOff

\def\Title{Central Extensions of Word Hyperbolic Groups}
\def\Authors{Walter D. Neumann and Lawrence Reeves}
\def\eval#1{\SOverline{#1}}
\def\len{\mop{len}}
\SetAuthorHead{\Authors}
\SetTitleHead{\Title}

\def\vale#1{\eval{#1}}
\def\valg#1{\pi(\eval{#1})}
\def\length#1{\mop{len}(#1)}
\def\area#1{\mop{area}(#1)}

\H \Title\par\Authors\par

Thurston has claimed (unpublished) that central extensions of
word hyperbolic groups by finitely generated abelian groups are
automatic. We show that they are in  fact biautomatic. Further, we show that
every 2-dimensional cohomology class on a word hyperbolic group can be
represented by a bounded 2-cocycle. This lends weight to the claim of Gromov
that for a word hyperbolic group, the cohomology in every dimension is
bounded.

Our results apply more generally to virtually central extensions.  We
build on the ideas presented in \cite{NR}, where the general problem
was reduced to the case of central extensions by $\bbZ$ and was solved
for Fuchsian groups. Some special cases of automaticity or biautomaticity
in this case had previously been proved in
\cite{ECHLPT},
\cite{Sha}, and \cite{Ge}.  

The new ingredient is a maximising
technique inspired by work of Epstein and Fujiwara. Beginning with an
arbitrary finite generating set for a central extension by $\bbZ$,
this maximising process is used to obtain a section which, in the
language of \cite{NR}, corresponds to a ``regular 2-cocycle'' on the
hyperbolic group G, and can be used to obtain a biautomatic structure
for the extension. Since central extensions correspond to
2-dimensional cohomology classes, it follows that every such class can
be represented by a regular 2-cocycle.  Using the geometric properties
of $G$, we then further modify this cocycle to obtain a bounded
representative for the original cohomology class.

We also discuss the relations between various concepts of ``weak
boundedness'' of a 2-cocycle on an arbitrary finitely generated group,
related to quasi-isometry properties of central extensions. For
cohomology classes these weak boundedness concepts are shown to be
equivalent to each other.  We do not know if a weakly bounded
cohomology class must be bounded.

\HH 1.~~Preliminaries

Let $G$ be a finitely generated group and $X$ a finite set which maps
to a monoid generating set of $G$. The map of $X$ to $G$ can be
extended in the obvious way to give a monoid homomorphism of $X^*$
onto $G$ which will be denoted by $w\mapsto\eval w$.  For convenience
of exposition we will always assume our generating sets are
\Em{symmetric}, that is, they satisfy $\eval X=\eval X^{-1}$.  If
$L\subset X^*$ then the pair consisting of $L$ and the evaluation map
$L\to G$ will be called a \Em{language on $G$}.  Abusing terminology,
we will often suppress the evaluation map and just call $L$ the
language on $G$ (but therefore, we may use two letters, say $L$ and
$L'$, to represent the same language $L\subset X^*$ with two different
evaluation maps to two different groups).  A language on $G$ is a
\Em{normal form} if it surjects to $G$.

A \Em{rational structure} for $G$ is a normal form $L\subset X^*$ for
$G$ which is a regular language (i.e., the set of accepted words for
some finite state automaton).  A subset $S$ of $G$ is then
\Em{$L$-rational} if the set of $w\in L$ with value in $S$ is a
regular language.

The \Em{Cayley graph} $\Gamma_X(G)$ is the directed graph with vertex
set $G$ and a directed edge from $g$ to $g\eval x$ for each $g\in G$
and $x\in X$; we give this edge a label $x$. 

Each word $w \in X^*$ defines a path $[0,\infty)\to\Gamma$ in the
Cayley graph $\Gamma=\Gamma_X(G)$ as follows (we denote this path also
by $w$): $w(t)$ is the value of the $t$-th initial segment of $w$ for
$t=0,\ldots,\len(w)$, is on the edge from $w(s)$ to $w(s+1)$ for
$s<t<s+1\le\len(w)$ and equals $\eval w$ for $t\ge\len(w)$. We refer
to the translate by $g\in G$ of a path $w$ by $gw$.

Let $\delta\in \Bbb N$. Two words $v,w\in X^*$ \Em{synchronously
$\delta$-fellow-travel} if the distance $d(w(t),v(t))$ never exceeds
$\delta$.  They \Em{asynchronously $\delta$-fellow-travel} if there
exist non-decreasing proper functions $t\mapsto t', t\mapsto t''\colon
[0,\infty)\to[0,\infty)$ 
such that $d({v(t')},{w(t'')})\le\delta$ for
all $t$.

A rational structure $L$ for $G$ is a \Em{synchronous} resp.\
\Em{asynchronous automatic structure} if there is a constant $\delta$
such that any two words $u,v\in L$ with $d(\eval u,\eval v)\le 1$
synchronously resp.\ asynchronously fellow-travel. A synchronous
automatic structure $L$ is \Em{synchronously biautomatic} if there is
a constant $\delta$ such that if $v,w\in L$ satisfy $\eval
w=\eval{xv}$ with $x\in X$ then $\eval xv$ and $w$ synchronously
$\delta$-fellow-travel.  See \cite{NS1} for a discussion of the
relationship of these definitions with those of \cite{ECHLPT}.  In
particular, as discussed there, if $L\to G$ is finite-to-one, then the
definitions are equivalent; by going to a sublanguage of $L$ this can
always be achieved.  

Two rational structures $L_1$ and $L_2$ on $G$ are \Em{equivalent} if
there is a $\delta$ such that each $w\in L_1$ asynchronously
$\delta$-fellow-travels some $v\in L_2$ with $\eval w=\eval v$ and vice
versa.  A subset of $G$ is then $L_1$-rational if and only if it is
$L_2$-rational (see \cite{NS2}).  A word-hyperbolic group has a unique
(bi)-automatic structure up to equivalence, given by the language of
geodesic words for any finite generating set (see \cite{NS1}).

\skipaline
If $G$ is a group and $A$ an abelian group then an extension 
$$0\to A\to E\to G\to 1$$
is called \Em{virtually central} if the induced map
$G\to\mop{Aut}(A)$ has finite image (we say the action of $G$ on
$A$ \Em{is finite}). In \cite{NR} the problem of
finding a biautomatic structure on $E$ is reduced to the case of
central extensions by $A=\bbZ$.  We describe this reduction now.

We write $A$ additively and
we denote the action of an element $g\in G$ on $A$ by $a\mapsto
a^g$.
Choose a section $s:G\to E$.  Then a general element of $E$ has
the form $s(g)\iota(a)$ with $g\in G$ and $a\in A$ and the group
structure in $E$ is given by a formula
$$s(g_1)\iota(a_1)s(g_2)\iota(a_2)
  =s(g_1g_2)\iota(a_1^{g_2}+a_2+\sigma(g_1,g_2)),$$
where $\sigma\colon G\times G\to A$ is a 2-cocycle on $G$ with
coefficients in $A$.  Changing the choice of section
changes the cocycle $\sigma$ by a coboundary.  Conversely, given a
cocycle $\sigma$, that is, a function $G\times G\to A$ which satisfies the
cocycle relation $$\sigma(g,h_1h_2)=\sigma(g,h_1)^{h_2} +
\sigma(gh_1,h_2)-\sigma(h_1,h_2),$$  
the above multiplication rule defines a virtually
central extension of $G$ by $A$.

\df Definition

Suppose $G$ has finite generating set $X$ and $L\subset X^*$ is an
asynchronous automatic structure on $G$.  We say a 2-cocycle $\sigma$
as above is \Em{weakly bounded} if 
 \Roster
 \ItemItem{1.} The sets $\sigma(g,G)$ and $\sigma(G,g)$ are finite for
each $g\in G$ (equivalently, $\sigma(X,G)$ and $\sigma(G,X)$ are both
finite --- this equivalence follows from the cocycle relation);
\endRoster
\noindent and is \Em{$L$-regular} if in addition
 \Roster
 \ItemItem{2.} For each $h\in G$ and $a\in A$ the
subset $\{g\in G: \sigma(g,h)=a\}$ is an $L$-rational subset of $G$
(it suffices to require this for $h\in X$ by \cite{NR, Lemma 2.1}).
\endRoster
\noindent A cohomology class in $H^2(G;A)$ is \Em{$L$-regular} if it
can be represented by an $L$-regular cocycle.  The term ``weakly
bounded'' reflects the standard terminology of \Em{bounded} for a
cocycle that satisfies $\sigma(G,G)$ finite.
\enddf

The following lemma is elementary:

\th Lemma 1.1

\Item{1.}$L$-regular cohomology classes form a subgroup
$H_L^2(G;A)\subset H^2(G;A)$ which only depends on the equivalence
class of $L$.
\Item{2.} If $A\to B$ is an equivariant map of finitely generated
abelian groups with finite $G$-actions then the induced map
$H^2(G;A)\to H^2(G;B)$ maps $H_L^2(G;A)$ to $H_L^2(G;B)$.\qed
\endth  

The following is proved in \cite{NR}.

\th Theorem 1.2

Let $0\to A\to E\to G\to 1$ be a virtually central extension with $A$
finitely generated.  Then $E$ has a biautomatic structure if and only
if $G$ has a biautomatic structure $L$ for which the cohomology class
of the extension is in $H^2_L(G;A)$.\qed
\endth

Our first main result can thus be formulated:

\th Theorem 1.3

Suppose $G$ is word hyperbolic and $L$ is the biautomatic structure on
$G$ (which is unique up to equivalence).  If $A$ is any finitely
generated abelian group with finite $G$-action then
$H^2_L(G;A)=H^2(G;A)$. In particular, any virtually central extension
of $G$ by $A$ is biautomatic.
\endth

In the remainder of this section we show that this theorem follows
from the special case that $A=\bbZ$ with trivial $G$-action. Assume
therefore that it is always true in this case.

Recall (e.g., \cite{NS1}) that if $H<G$ is a subgroup of finite index
and $L$ is a biautomatic structure on $G$ then there is an induced
biautomatic structure $L_H$ on $H$ which is unique up to equivalence.

\th Lemma 1.4. {\rm(\cite{NR, Corollaries 2.5 and 2.7})}

Let $L$ be a biautomatic structure on $G$.  Suppose $\alpha\in
H^2(G;A)$.

\Item{1.}Let $H<G$ be a subgroup of finite index and suppose $\beta\in
H^2(H;A)$ is the restriction of $\alpha$ to $H$. Then $\alpha$ is
$L$-regular if and only of $\beta$ is $L_H$-regular.  
\Item{2.}  A
non-zero multiple of $\alpha$ is $L$-regular if and only if $\alpha$
is.  \qed
\endth

We restrict to the kernel $H$ of the $G$-action on $A$. For $\alpha\in
H^2(G;A)$ we denote its restriction by $\beta\in H^2(H;A)$. By part 1
of the lemma it suffice to show $\beta$ is regular.  Write $A$ as a
direct sum of a finite group $F$ and copies of $\bbZ$ as follows:
$A=F\oplus\coprod_{i=1}^n \bbZ$. Then
$H^2(H;A)=H^2(H;F)\oplus\coprod_{i=1}^n H^2(H;\bbZ)$.  For each
$j=1,\ldots,n$ we can form $\bbZ_j=A/(F\oplus\coprod_{i\ne j}\bbZ)$
and the induced class $\beta_j\in H^2(H;\bbZ_j)$ is regular by
assumption.  We can embed $\bbZ_j$ in $A$ as the $j$-th summand $\bbZ$
and thus lift $\beta_j$ to a regular class $\beta'_j\in H^2(H;A)$.
Now $\beta$ differs from the regular class
$\beta':=\sum_{i=1}^n\beta'_j$ by torsion, so $\beta$ and $\beta'$
have a common nontrivial multiple. Thus $\beta$ is regular by part 2
of the Lemma.\qed

\HH 2.~~Central Extensions of Hyperbolic Groups are Biautomatic

In view of the preceding discussion we only need show that it $G$ is a
hyperbolic group and
$$0\to \bbZ\mathop{\to}\limits^\iota E\mathop{\to}\limits^\pi G\to 1$$
is a central extension of $G$ then $E$ is biautomatic.  

Let $X$ be a finite set which maps to a symmetric generating set for
$E$.  Then $X$ also maps to a generating set for $G$. We use the
notations $\vale x$ and $\valg x$ for evaluation of $x$ into $E$ and
$G$.  Fix some finite presentation $\langle X | R \rangle$ for $G$.

For each $g\in G$ the fibre $\pi^{-1}(g)$ has a total ordering
determined by the usual ordering on $\bbZ$.  It therefore makes sense
to talk about existence of a ``maximum'' for a subset of this fiber.

\th Lemma 2.1
 
There exists $C>0$ such that for all $g\in G$
$$\max\{\vale w \iota(-C\length w):w\in X^*,\valg{w}=g\}$$
exists. Moreover, there then exists $\lambda>0$ so that any word $w$
which achieves the maximum defines a $(\lambda, 0)$-quasigeodesic in
$G$.
\endth

\pf Proof

Let $T\in\bbZ$ be given by $\iota(T)=\max\{\vale{r}: r^{\pm1}\text{ is a
relator in }G\}$, and let $K$ be the constant
in the linear isoperimetric inequality for $G$, that is, if $\valg w=
1$ then the combinatorial area bounded by $w$ is less than $K\length
w$. In particular, we then have $\vale w\le
\iota(T\area{w})\le \iota(TK\length w)$. Choose $C\in\bbZ$
with $C>TK$.

Fix an element $g\in G$, and suppose that $w,\gamma\in X^*$ satisfy
$\valg w=\valg \gamma=g$. We have
$$\eqalignno{\vale \gamma ^{-1}\vale w 
&\le \iota(TK\length{\gamma^{-1}w})\cr
&= \iota(TK\length \gamma)\iota(TK\length w)&(1)\cr}$$
So,
$$\eqalign{\vale w \iota(-C\length w)&\le
\vale w \iota(-TK\length w)\cr  
&\le \vale \gamma\iota(TK\length{\gamma})\cr}$$
and it follows that the maximum exists in the lemma.  
We shall call a word which achieves the maximum in Lemma 2.1 a
\Em{maximising word}. 

Now assume that $w$ is maximising and that $\gamma$ defines a geodesic
in $G$. Then
$$\vale \gamma \iota(-C\length
\gamma)\le \vale w \iota(-C\length w)$$
so 
$$\iota(C(\length w - \length \gamma)) \le \vale \gamma^{-1}\vale w \eqno(2)$$
Combining $(1)$ and $(2)$ gives $C(\length w - \length \gamma)\le
TK(\length \gamma-\length w)$, and hence
$$\length w \le ((C+TK)/(C-TK))\length \gamma.$$  
Since a subword of a maximising word is a itself a maximising word,
this says that $w$ is $(\lambda,0)$-quasigeodesic in $G$
with $\lambda=((C+TK)/(C-TK)$.\qed

Define a new section $\rho\colon G\to E$ 
by setting $\rho(g)=\vale w\iota(-C\length w )$, where $w$ is a
maximising word such that $\valg w=g$.  The following is the main
result of this section.

\th Proposition 2.2

The cocycle  defined by the section $\rho$ is a regular cocycle.
\endth

\pf Proof

For elements $h_1,h_2$ in a common fibre of $E\to G$ we shall denote
by $h_1-h_2$ the integer with $\iota(h_1-h_2)=h_1h_2^{-1}$.  The
cocycle is then defined by the formula
$$\sigma(g_1,g_2)=\rho(g_1)\rho(g_2)-\rho(g_1g_2).$$
The weak boundedness of $\sigma$ is therefore the statement that
$\rho(g)\rho(x)-\rho(gx)$ and $\rho(x)\rho(g)-\rho(xg)$ are bounded
functions of $(g,x)\in G\times X$.  The following lemma clearly
implies this, since $X$ is finite.

\th Lemma 2.3

Suppose $g\in G$ and $x\in X$.  Then
\Item{1.} 
$|\rho(g)\vale x-\rho(gx)|\le C$.
\Item{2.}$|\vale x\rho(g)-\rho(xg)|\le C$.
\endth

\pf Proof

We prove part 1. Choose maximising words $w_1$ and $w_2$ with $\valg
{w_1}=g$ and $\valg {w_2}=gx$. Since ${w_2}$ is maximising and
$\valg{w_1x}=\valg{w_2}$ we have
$$
\vale{w_1}\iota(-C\length {w_1} )\vale x \iota(-C)\iota(h)= \vale{w_2}
\iota(-C\length {w_2}),
$$ 
for some $h\ge 0$. By symmetry, since $\valg{w_1}=\valg{w_2x^{-1}}$,
we have
$$
\vale {w_2}\iota(-C\length {w_2} )\vale x^{-1} \iota(-C)\iota(k)=
\vale {w_1} \iota(-C\length {w_1}),
$$ 
for some $k\ge 0$. Comparing, we see $-C\le h-C=C-k\le C$
and the result follows. Part 2 is completely analogous.\qed

Let $L\subset X^*$ be the language of maximising words.  Since $L$
consists of $(\lambda,0)$ quasigeodesics in $G$ and quasigeodesics
fellow-travel geodesics (see \cite{G} or \cite{Sho}), the language $L$
has the asynchronous fellow traveller property.  We shall show next
that it is a regular language, so it is an asynchronous automatic
structure which is equivalent to the geodesic automatic structure on
$G$. (It is not necessarily itself synchronously automatic --- for
example it fails to be synchronous if $G$ is a surface group and $E$ a
non-trivial extension.)  Let $\delta$ be the fellow-traveller constant
for $L$.

We first note that if $u\in X^*$ is a word that is not in $L$ then
this fact can be recognised via a $\delta$-fellow-traveller.  That is,
there is a $v\in X^*$ with $\valg u= \valg v$ such that $v$ and $u$
$\delta$-fellow-travel and $v$ is ``better'' than $u$ in that $\vale
u\iota(-C\length u)< \vale v\iota(-C\length v)$.  Indeed, we obtain
$v$ by taking the shortest initial segment $u_1$ of $u$ that fails to
be in $L$ (i.e., is not maximising) and replacing it by a maximising
word $v_1$.  This $u_1$ and $v_1$ $\delta$-fellow travel because $u_1$
has the form $u_0x$ with $u_0\in L$ and $x\in X$.

Given this observation, called ``falsification by fellow-traveller''
in \cite{NS2}, the regularity of $L$ is a standard argument,
originally due to Cannon (see \cite{NS2, Proposition 4.1}).  We sketch
it here for completeness.

To test if a path $u$ is maximising, as we move along $u$ we must keep
track at each time $t$ of what points $u(t)h$ in a
$\delta$-neighborhood of $u(t)$ have been reached by paths $v$ that
asynchronously $\delta$-fellow travel $u$.  Moreover, for each such
$g$ we must record the difference in ``height'' of $\vale
{u(t)}\iota(-Ct)$ and $\vale{ vw}\iota(-C \length {vw})$ where $w$ is
a geodesic with value $h^{-1}$. If this difference is larger than
$2\delta C$ then $v$ is not a viable competitor to $u$ and we can
ignore it (following $u$ to $u(t)$ and then $w^{-1}$ to $\valg v$ does
better than $v$).  If this difference is negative then the path $vw$
does better than $u$, so $u$ is not in $L$.  We are thus only
interested in differences that lie in the interval $[0,2\delta C]$.
At each point $h$ of the ball $B(\delta)$ of radius $\delta$ about $1$
in the Cayley graph we record the ``best'' value of this difference
for paths $v$ ending at $u(t)h$.  Thus the information that we must
keep track of is a function $\phi\colon
B(\delta)\to\{0,1,\ldots,2\delta C\}$.  We can build a finite state
automaton with the set of such maps as states plus one ``fail state,''
since, as one moves along $u$, the state at any point of the path $u$
can clearly be computed from its value at the previous point and the
last letter of $u$ read.  Once we reach the end of $u$, we know $u$ is
in $L$ if and only if we have not yet reached the fail state.  Thus
$L$ is recognised by a finite state automaton so it is a regular
language.

We now lift $L$ via the section $\rho$ as follows.  For $x\in X$
denote $x':=\vale x\iota(-C)$.  Denote by $L'$ the language $L$ with
the evaluation map $x\mapsto x'$ instead of $x\mapsto \vale x$.  If
$v=x_1x_2\ldots x_n\in L'$ then the initial segments of $v$ have
values $\rho(\valg{x_1})$, $\rho(\valg{x_1x_2}),\ldots,$
$\rho(\valg{x_1\ldots x_n})$.  $L'$ is a regular language mapping onto
the image of our section $\rho$ with asynchronous fellow-traveller
property. As in the proof of Proposition 2.2 of \cite{NR}, this
implies that the cocycle $\sigma$ defined by $\rho$ is regular.
Indeed, for any $x\in X$ and $a\in \bbZ$ the set
$\{(w_1,w_2):w_1,w_2\in L,\ \rho(\vale{w_1})\rho(\vale x)\iota(-a)=
\rho(\vale{w_2})\}$ is the language of an asynchronous two-tape
automaton so its projection to the first factor is a regular language.
But that projection is $\{g:\sigma(g,x)=a\}$, so this set is rational,
as required.\qed

We have thus proved the biautomaticity of $E$.  It is worth recalling
from \cite{NR} that a biautomatic structure is found by taking a
finite-to-one automatic structure on $G$ (e.g., the geodesic language)
and lifting it using the section $\sigma$ to get a language $M$
evaluating onto $\sigma$ (the alphabet for $M$ consists of elements of
the form $\vale x\iota(-\sigma(g,x))$ with $g\in G$ and $x\in X$).
Then $M(\{t\}^*\cup\{t^{-1}\}^*)$ is the desired biautomatic
structure, where $t$ is a generator of $\bbZ$.  The language $L'$ of
lifted maximising words cannot be used instead of $M$ since it is not
generally synchronous.

\HH 3.~~ Bounded Cohomology

The main result of this section is

\th Theorem 3.1

Suppose $G$ is word hyperbolic and  $A$ is any finitely
generated abelian group with finite $G$-action then any class in 
$H^2(G;A)$ can be represented by a bounded  cocycle.
\endth

As in section 1, we can reduce to the case that $A=\bbZ$ with trivial
$G$-action.  Given a 2-cocycle, and hence some central extension, in
the previous section we defined a section $\rho$ determining a
cohomologous cocycle which was weakly bounded.
In this section we modify $\rho$ to define a section for which
the corresponding cocycle is bounded.  That is, given a weakly bounded
2-cocycle on a word hyperbolic group we produce a cohomologous bounded
2-cocycle.

\th Lemma 3.2

Suppose that $k$ is a fixed constant.  Then there is a constant $K(k)$
such that if $g,h\in G$ are such that $g$ lies within distance $k$ of
a geodesic path from 1 to $gh$, then $d_E(\rho(gh),\rho(g)\rho(h))\le
K(k)$.
\endth

\pf Proof

Our section $\rho$ has the property that if $g$ lies on the maximising
path from $1$ to $gh$ then $\rho(gh)=\rho(g)\rho(h)$.  The weak
boundedness of the cocycle associated to our section thus implies that
if $g$ lies a bounded distance from the maximising path to $gh$ then
$d_E(\rho(gh),\rho(g)\rho(h))$ is bounded.  Since maximising paths
fellow travel geodesics, the lemma follows.
\qed

To show boundedness of our cocycle $\sigma$ we shall need that the
section is symmetric. To ensure this we define a new section $q$ by
taking the average of $\rho(g)$ and $\rho(g^{-1})^{-1}$, that is,
$$q(g)=1/2(\rho(g)+\rho(g^{-1})^{-1}).$$
This section is defined in the pushout extension
$$0\to\bbQ\to E_\bbQ\to G\to 1$$
determined by the embedding $\bbZ\to\bbQ$.  Note that although the
fibres of $E_\bbQ\to G$ only have an affine structure, the average of
two points in a fiber is nevertheless well defined.  We shall show
that the cocycle defined by this section is bounded.  It
then follows easily that the cocycle defined by the section $[q]\colon
G\to E$ (integral part of $q$) is bounded.

It is clear from the definition that 

\th Lemma 3.3

For all $g\in G$, $q(g^{-1})=q(g)^{-1}$.
\qed
\endth

\th Proposition 3.4

The cocycle determined by $q$ is bounded.
\endth

\pf Proof

We need to show that there  is a bound on the values attained by
by $q(gh)^{-1}q(g)q(h)$, for $g,h\in G$. For $g\in G$ denote by $w_g$
a geodesic path from $1$ to $g$.  Consider the geodesic triangle in
$G$ with sides labelled by $w_{gh}$, $w_g$, and $gw_h$. Since $G$ is
hyperbolic there is a constant $k$ such that for all such triangles
there is a $c\in G$ such that $c$ lies within distance $k$ of each of
the sides.  Notice that $q$ inherits from $\rho$ the property
described in Lemma 3.2.  That is replacing $q(g)$ by $q(c)q(c^{-1}g)$
alters the value by a bounded amount. Similarly we replace $q(h)$ by
$q(g^{-1}c)q(c^{-1}gh)$, and $q(gh)$ by $q(c)q(c^{-1}gh)$. So,
$q(gh)^{-1}q(g)q(h)$ differs from
$(q(c)q(c^{-1}gh))^{-1}q(c)q(c^{-1}g)q(g^{-1}c)q(c^{-1}gh)=1$ by a
bounded amount.  \qed

\rk Remark

It is not hard to verify that the averaging process described in this
section to make our cocycle bounded does not destroy regularity, so,
in fact, any class in $H^2(G;A)$ has a bounded regular representative.
\endrk

\HH 4. Weakly bounded cocycles

In this section we do not assume $G$ is word hyperbolic.  We do assume
$G$ is finitely generated.  We shall call a 2-cocycle $\sigma$ on $G$
\Em{left weakly bounded} if $\sigma(g,G)$ is bounded for all $g\in G$
and \Em{right weakly bounded} if $\sigma(G,g)$ is bounded for all
$g\in G$.  So a cocyle is weakly bounded if and only if it is both
left and right weakly bounded.

In \cite{NR} we pointed out that if a central extension $0\to A\to
E\to G\to 1$ is given by a right weakly bounded cocycle then $E$ is
quasi-isometric to $G\times A$ . In fact, the corresponding section
$s\colon G\to E$ is quasi-isometric and induces a quasi-isometry
$G\times A\to E$ by $(g,a)\mapsto s(g)\iota(a)$.  Conversely, if one
has a quasi-isometric section then the corresponding cocycle is right
weakly bounded.  We asked there about the relation between the
concepts right weakly bounded, weakly bounded, and bounded.  We still
do not know if a weakly bounded cohomology class is always bounded,
but we have:

\th Theorem 4.1

A cocycle which is either right or left weakly bounded is always
cohomologous to one which is both right and left weakly bounded.
\endth

\pf Proof

The same argument as in section 1 can be used to reduce to the case
that $A=\bbZ$.  We assume therefore that we have a central extension
$0\to \bbZ\to E\to G\to 1$ given by a right weakly bounded cocycle
(the left weakly bounded case is completely analogous).  As usual, we
assume we have a symmetric generating set $X$ for $E$, and we also
consider it as a generating set for $G$. Let $s:G\to E$ be a section
whose cocycle $\sigma$ is right weakly bounded.  Since changing
$\sigma$ at finitely many points does not destroy right weak
boundedness, we may assume $s(1)=1$ and $s(\valg x)=\vale x$ for $x\in
X$.  Let $K$ be a bound with $|\sigma(g,\valg x)|\le K$ for $g\in G$
and $x\in X$.  Then if $w=x_1x_2\ldots x_n$ is any word we have
$s(\pi(w(t)))=s(\pi(w(t-1)))\vale{x_t}\iota(-\sigma(w(t-1),\valg{x_t}))$
for each $t$.  Thus if $w$ represents $1$ in $G$ we have by induction
$1=s(\valg
w)=\vale{x_1x_2...x_n}\iota(-\sum_{t=1}^n\sigma(w(t-1),\valg{x_t}))$.
Thus $\vale w=\vale{x_1\ldots x_n}$ is bounded by $nK$ in absolute
value. Thus if we choose $C>K$ the proof of Lemma 2.1 goes through.
The section defined by maximising paths gives, as in section 2, a
weakly bounded cocycle.\qed


\Bib        Bibliography

\rf {ECHLPT} D.~B.~A. Epstein, J.~W. Cannon, D.~F. Holt, S.~V.~F. Levy,
M.~S. Paterson, and  W.~P. Thurston, Word processing in groups, Jones and
Bartlett, 1992.
\rf {Ge} S.~M.~Gersten, Bounded cocycles and combings of groups, 
Int. J. of Algebra and  Computation {\bf 2} (1992), no. 3, 307--326.

%

\rf {G} M. Gromov, Hyperbolic groups, 1987, 75--263, Essays in Group
Theory, Springer Verlag, M.S.R.I.  series Vol 8, S. M. Gersten,
editor.

%
\rf {NR} W.~D.~Neumann and L.~Reeves, Regular cocycles and biautomatic
structures, Int. J. Alg. and Comp. (to appear).

\rf {NS1} W.~D.~Neumann and M.~Shapiro, Equivalent automatic structures and
their boundaries, Int. J. of Algebra and  Computation {\bf 2} (1992), no. 4,
443--469.

\rf {NS2} W.~D.~Neumann and M.~Shapiro, Automatic structures, rational growth,
and geometrically finite groups. Invent. Math. (to appear).

\rf{Sha} M.~Shapiro, Automatic structures and graphs of groups, in:
``Topology `90, Proceedings of the Research Semester in Low
Dimensional Topology at Ohio State,'' Walter de Gruyter Verlag,
Berlin - New York (1992), 355--380.

\rf{Sho} H. Short (editor), Notes on word hyperbolic groups,
in ``Group theory from a geometric viewpoint, conference held
in I.C.T.P. Trieste, March 1990'', edited by E. Ghys, A. Haefliger
and A. Verjovsky, (World Scientific, 1993)

\endBib

\bye